\documentstyle[amscd,amssymb,verbatim,epsf]{amsart}
\begin{document}

\theoremstyle{plain}
\newtheorem{Thm}{Theorem}
\newtheorem{Cor}{Corollary}
\newtheorem{Con}{Conjecture}
\newtheorem{Main}{Main Theorem}
\newtheorem{Exmp}{Example}
\newtheorem{Lem}{Lemma}
\newtheorem{Prop}{Proposition}

\theoremstyle{definition}
\newtheorem{Def}{Definition}
\newtheorem{Note}{Note}

\theoremstyle{remark}
\newtheorem{notation}{Notation}
\renewcommand{\thenotation}{}

\errorcontextlines=0
\numberwithin{equation}{section}
\renewcommand{\rm}{\normalshape}%

\title[Upward and downward statistical continuities]%
   {Upward and downward statistical continuities}
\author{Huseyin Cakalli\\
          Maltepe University, Maltepe, Istanbul-Turkey}
\address{Huseyin Cakalli\
          Maltepe University, Department of Mathematics, Marmara E\u{g}\.{I}t\.{I}m K\"oy\"u, TR 34857, Maltepe, \.{I}stanbul-Turkey Phone:(+90216)6261050 ext:2248, fax:(+90216)6261113
}
\email{hcakalli@@maltepe.edu.tr; hcakalli@@gmail.com}

\date{\today}

\begin{abstract}
A real valued  function $f$ defined on a subset $E$ of $\textbf{R}$, the set of real numbers, is statistically upward continuous if it preserves statistically upward half quasi-Cauchy sequences, is statistically downward continuous if it preserves statistically downward half quasi-Cauchy sequences; and a subset $E$ of $\textbf{R}$, is statistically upward compact if any sequence of points in $E$ has a statistically upward half quasi-Cauchy subsequence, is statistically downward compact if any sequence of points in $E$ has a statistically downward half quasi-Cauchy subsequence where a sequence $(x_{n})$ of points in $\textbf{R}$ is called statistically upward half quasi-Cauchy if
\[
\lim_{n\rightarrow\infty}\frac{1}{n}|\{k\leq n: x_{k}-x_{k+1}\geq \varepsilon\}|=0
\]
is statistically downward half quasi-Cauchy if
\[
\lim_{n\rightarrow\infty}\frac{1}{n}|\{k\leq n: x_{k+1}-x_{k}\geq \varepsilon\}|=0
\]
for every $\varepsilon>0$.
We investigate statistically upward continuity, statistically downward continuity, statistically upward half compactness, statistically downward half compactness and prove interesting theorems. It turns out that uniform limit of a sequence of statistically upward continuous functions is statistically upward continuous, and uniform limit of a sequence of statistically downward continuous functions is statistically downward continuous.

\end{abstract}

\maketitle

\section{Introduction}
\normalfont
 Connor and Grosse-Erdmann \cite{ConnorandGrosse-ErdmannSequentialdefinitionsofcontinuityforrealfunctions} gave sequential definitions of continuity for real functions calling $G$-continuity instead of $E$-continuity and their results cover the earlier works related to $E$-continuity where a method of sequential convergence, or briefly a method, is a linear function $G$ defined on a linear subspace of $s$, denoted by $c_{G}$, into $\textbf{R}$. A sequence $\textbf{x}=(x_{n})$ is said to be $G$-convergent to $\ell$ if $\textbf{x}\in c_{G}$ and $G(\textbf{x})=\ell$. In particular, $\lim$ denotes the limit function $\lim \textbf{x}=\lim_{n}x_{n}$ on the linear space $c$. A function $f$ is called $G$-continuous at a point $u$ provided that whenever a sequence $\textbf{x}=(x_{n})$ of terms in the domain of $f$ is $G$-convergent to $u$, then the sequence $f(\textbf{x})=(f(x_{n}))$ is $G$-convergent to $f(u)$. A method $G$ is called regular if every convergent sequence $\textbf{x}=(x_{n})$ is $G$-convergent with $G(\textbf{x})=\lim \textbf{x}$. A method is called subsequential if whenever $\textbf{x}$ is $G$-convergent with $G(\textbf{x})=\ell$, then there is a subsequence $(x_{n_{k}})$ of $\textbf{x}$ with $\lim_{k} x_{n_{k}}=\ell$. Recently, Cakalli  gave new sequential definitions of compactness and slowly oscillating compactness in \cite{CakalliSequentialdefinitionsofcompactness}, and \cite{CakalliSlowlyoscillatingcontinuity}, respectively. The notion of $N_\theta$ convergence was introduced, and studied by Freedman, Sember, and M. Raphael in \cite{FreedmanandSemberandRaphaelSomecesarotypesummabilityspaces}. Using the idea of Sember and Raphael, Fridy and Orhan introduced lacunary statistical convergence (\cite{FridyandOrhanLacunarystatisconvergence} and  \cite{CakalliLacunarystatisticalconvergenceintopgroups}).

A subset $E$ of $\textbf{R}$, the set of real numbers, is compact if and only if any sequence of points in $E$ has a convergent subsequence whose limit in $E$. A subset $E$ of $\textbf{R}$ is bounded if and only if any sequence of points in $E$ has a Cauchy subsequence. Boundedness coincides with not only ward compactness (\cite{CakalliForwardcontinuity}), but also either of the following kinds of compactnesses, slowly oscillating compactness (\cite{CakalliSlowlyoscillatingcontinuity}), statistical ward compactness (\cite[Lemma 2]{CakalliStatisticalwardcontinuity}), lacunary statistical ward compactness, (\cite[Theorem 3]{CakalliStatisticalquasiCauchysequences}) $N_{\theta}$-ward compactness (\cite[Theorem 3.3]{CakalliNthetawardcontinuity}).
Two of our results in this paper state necessary and sufficient conditions for below boundedness and above boundedness of a subset of $\textbf{R}$.
Using the idea of continuity of a real function in terms of sequences, many kinds of continuities were introduced and investigated, not all but some of them we recall in the following: slowly oscillating continuity  (\cite{CakalliSlowlyoscillatingcontinuity}), quasi-slowly oscillating continuity (\cite{DikandCanak}), ward continuity (\cite{CakalliForwardcontinuity}, \cite{BurtonandColemanQuasiCauchysequences}), $\delta$-ward continuity (\cite{CakalliDeltaquasiCauchysequences}), statistical ward continuity,  (\cite{CakalliStatisticalwardcontinuity}), and $N_{\theta}$-ward continuity (\cite{CakalliNthetawardcontinuity}, \cite{CakalliandKaplanAstudyonNthetaquasiCauchysequences}) which enabled some authors to obtain conditions on the domain of a function for some characterizations of uniform continuity in terms of sequences in the sense that a function preserves a certain kind of sequences (see \cite{Vallin} and \cite{BurtonandColemanQuasiCauchysequences}).

The purpose of this paper is to introduce the concepts of statistically upward and statistically downward continuities, and prove interesting theorems.

\maketitle

\section{Preliminaries}

\normalfont{}
Throughout the paper, $s$ and $c$ will denote the set of all sequences, and the set of convergent sequences of points in $\textbf{R}$.

The concept of statistical convergence is a generalization of the usual notion of convergence that, for real-valued sequences, parallels the usual theory of convergence. A sequence $(x_{k})$ of points in \textbf{R} is called statistically convergent to an element $\ell$ of \textbf{R}  if for each
$\varepsilon$
\[
\lim_{n\rightarrow\infty}\frac{1}{n}|\{k\leq n: |x_{k}-\ell|\geq{\varepsilon}\}|=0,
\] and this is denoted by $st-\lim_{k\rightarrow\infty}x_{k}=\ell$ (see \cite{FridyOnstatisticalconvergence}, \cite{MaioKocinac}, \cite{CakalliAstudyonstatisticalconvergence}, \cite{CasertaandKocinacOnstatisticalexhaustiveness},  \cite{CasertaMaioKocinacStatisticalConvergenceinFunctionSpaces}, \cite{CakalliandKhan}, \cite{CakalliandKhanSummabilityintopologicalspaces}  and \cite{MursaleenLamdastatisticalconvergence}). This defines a regular subsequential method of sequential convergence, i.e. $G(\boldsymbol{\alpha}):= st-\lim_{k\rightarrow\infty}x_{k}$.

A lacunary sequence $\theta=(k_{r})$ is an increasing sequence $\theta=(k_{r})$ of positive integers such that $k_{0}=0$ and $h_{r}:k_{r}-k_{r-1}\rightarrow\infty$. The intervals determined by $\theta$ will be denoted by $I_{r}=(k_{r-1}, k_{r}]$, and the ratio $\frac{k_{r}}{k_{r-1}}$ will be abbreviated by $q_{r}$. Sums of the form $\sum^{k_{r}}_{k_{r-1}+1}|x_{k}|$ frequently occur, and will often be written for convenience as $\sum^{}_{k\in{I_{r}}}|x_{k}|$. Throughout this paper, we will assume that $\lim inf_{r}\; q_{r}>1$.

A sequence $(x_{k})$ of points in \textbf{R} is called $N_\theta$-convergent to an element $\ell$ of \textbf{R} if
\[
\lim_{r\rightarrow\infty}\frac{1}{h_{r}}\sum^{}_{k\in{I_{r}}}|x_{k}-\ell|=0,
\]
and it is denoted by $N_{\theta}-lim\; x_{k}=\ell$. This defines a regular method of sequential convergence, i.e. $G(\boldsymbol{\alpha}):=N_{\theta}-lim\; x_{k}$. Any convergent sequence is $N_{\theta}$-convergent, but the converse is not always true.

A real sequence $(x_{k})$ is called lacunary statistically convergent to an element $\ell$ of $\textbf{R}$ if
\[
\lim_{r\rightarrow\infty}\frac{1}{h_{r}} |\{k\in I_{r}: |x_{k}-\ell|\geq{\epsilon} \}|=0,
\]
for every $\epsilon>0$ where $\theta=(k_{r})$ is a lacunary sequence. This defines a regular subsequential method of sequential convergence, i.e. $G(\boldsymbol{\alpha}):=S_{\theta}-lim\; x_{k}$. Any convergent sequence is $S_{\theta}$-convergent, but the converse is not always true. For example, limit of the sequence of the ratios of Fibonacci numbers converge to the golden mean. This ensures the regularity of lacunary sequential method obtained via the sequence of Fibonacci numbers, i.e. $\theta=(k_{r})$ is the lacunary sequence defined by writing $k_{0}=0$ and $k_{r}=F_{r+2}$ where $(F_{r})$ is the Fibonacci sequence, i.e. $F_{1}=1$, $F_{2}=1$, $F_{r}= F_{r-1} + F_{r-2}$ for $r\geq 3$ (\cite{CakalliNthetawardcontinuity}).

The concept of a Cauchy sequence involves far more than that the distance between successive terms is tending to zero. Nevertheless, sequences which satisfy this weaker property are interesting in their own right.  A sequence $(x_{n})$ of points in \textbf{R} is quasi-Cauchy if $(\Delta x_{n})$ is a null sequence where $\Delta x_{n}=x_{n}-x_{n+1}$. These sequences were named as quasi-Cauchy by Burton and Coleman \cite[page 328]{BurtonandColemanQuasiCauchysequences}, while they  were called as forward convergent to $0$ sequences in \cite[page 226]{CakalliForwardcontinuity}.
A subset of \textbf{R} is compact if and only if it is closed and bounded. A subset $E$ of  \textbf{R} is bounded if $|a|\leq{M}$ for all $a \in{A}$  where  $M$ is a positive  real constant number. This is equivalent to the statement that any sequence of points in $E$ has a quasi-Cauchy subsequence.
A sequence $(x_{n})$ of points in \textbf{R} is slowly oscillating if
$$
\lim_{\lambda \rightarrow 1^{+}}\overline{\lim}_{n}\max _{n+1\leq
k\leq [\lambda n]} |
  x_{k}  -x_{n} | =0
$$ where $[\lambda n]$ denotes the integer part of $\lambda n$ (see \cite[Definition 2 page 947]{DikFandDikMandCanakIApplicationsofsubsequentialTauberiantheorytoclassicalTauberiantheory}).
Any subsequence of a Cauchy sequence is Cauchy. The analogous property fails for quasi-Cauchy sequences, and fails for slowly oscillating sequences as well. A counterexample for the case, quasi-Cauchy, is the sequence $(a_n)=(\sqrt{n})$ with the subsequence $(a_{n^{2}})=(n)$. A counterexample for the case slowly oscillating is the sequence $(log_{10} n)$ with the subsequence (n). Furthermore we give more examples without neglecting: the sequences $(\sum^{\infty}_{k=1}\frac{1}{n})$, $(ln\;n)$, $(ln\;(ln\;n))$, $(ln\;( ln\;( ln n)))$,..., $(ln\;( ln\;( ln (...( ln n)...))$   and combinations like that are all slowly oscillating, but not Cauchy. The bounded sequence $(cos (6 log(n + 1)))$ is slowly oscillating, but not Cauchy. The sequences $(cos(\pi \sqrt{n}))$ and $(\sum^{k=n}_{k=1}(\frac{1}{k})(\sum^{j=k}_{j=1}\frac{1}{j}))$ are quasi-Cauchy, but slowly oscillating(see \cite{CakalliSlowlyoscillatingcontinuity}, \cite{Vallin}, and   \cite{CakalliandSonmezSlowlyOscillatingContinuityinAbstractMetricSpaces}).

Now we recall the concepts of ward compactness, and slowly oscillating compactness: a subset $E$ of  \textbf{R} is called ward compact if any sequence of points in $E$ has a quasi-Cauchy subsequence (\cite{CakalliForwardcontinuity}). A subset $E$ of  \textbf{R} is called slowly oscillating compact if any sequence of points in $E$ has a slowly oscillating subsequence (\cite{CakalliSlowlyoscillatingcontinuity}).

\c{C}akall\i\; (\cite[page 594]{CakalliSequentialdefinitionsofcompactness}, \cite{CakalliOnGcontinuity}) gave a sequential definition of compactness, which is a generalization of ordinary sequential compactness, as in the following: a subset $E$ of \textbf{R} is $G$-sequentially compact if for any sequence $(x_{k})$ of points in $E$ there exists a subsequence $\textbf z$ of the sequence such that $G(\textbf{z})\in{A}$. His idea enables us obtaining new kinds of compactness via most of the non-matrix sequential convergence methods, as well as all matrix sequential convergence methods.

Palladino (\cite{PalladinoOnhalfcauchysequences}) introduced a concept of upward half Cauchness, and a concept of downward half Cauchyness as in the following:a sequence $(x_{n})$ of points in $\textbf{R}$ is called upward half Cauchy if for every $\varepsilon>0$ there exists an $n_{0}\in{\textbf{N}}$ so that $x_{n}-x_{m} <\varepsilon$ for $m \geq n \geq n_0$, downward half Cauchy if for every $\varepsilon>0$ there exists an $n_{0}\in{\textbf{N}}$ so that $x_{m}-x_{n} <\varepsilon$ for $m \geq n \geq n_0$. A sequence $(x_{n})$ is called half Cauchy if the sequence is either upward half Cauchy, or downward half  Cauchy, or both. It is clear that a sequence $(x_{n})$ is Cauchy if and only if it is both upward half Cauchy and downward half Cauchy.

Quasi-Cauchy sequences arise in diverse situations, and it is often difficult to determine whether or not they converge, and if so, to which limit. It is easy to construct a zero-one sequence such that the quasi-Cauchy average sequence does not converge. The usual constructions have a somewhat artificial feeling. Nevertheless, there are sequences which seem natural, have the quasi-Cauchy property, and do not converge. Using the idea of the definition of an upward half Cauchy sequence, a concept of upward half quasi-Cauchy sequence is introduced (see \cite{CakalliHalfquasi-Cauchysequences}, and \cite{CakalliUpwardanddownwardcontinuities}). A sequence $(x_{n})$ of points in $\textbf{R}$ is called upward half quasi-Cauchy if for every $\varepsilon>0$ there exists an $n_{0}\in{\textbf{N}}$ such that $x_{n}-x_{n+1} <\varepsilon$ for $n \geq n_0$.
A sequence $\textbf{x}=(x_{n})$ of points in $\textbf{R}$ is upward half Cauchy if and only if every subsequence of $\textbf{x}$ is upward half quasi-Cauchy.
A subset $E$ of $\textbf{R}$ is called upward compact if any sequence of points in $E$ has an upward half quasi-Cauchy subsequence.
A subset of $\textbf{R}$ is upward compact if and only if it is bounded below.
A subset $E$ of $\textbf{R}$ is called up half compact if any sequence of points in $E$ has an up half Cauchy subsequence (see \cite{CakalliHalfquasi-Cauchysequences}, and \cite{CakalliUpanddowncontinuities}).
A subset of $\textbf{R}$ is up half compact if and only if it is upward half compact.
A sequence $(x_{n})$ of points in $\textbf{R}$ is called downward half quasi-Cauchy if for every $\varepsilon>0$ there exists an $n_{0}\in{\textbf{N}}$ such that $x_{n+1}-x_{n} <\varepsilon$ for $n \geq n_0$ (\cite{CakalliHalfquasi-Cauchysequences}).
A sequence $(x_{n})$ of points in $\textbf{R}$ is downward half Cauchy if and only if every subsequence of $(x_{n})$ is downward half quasi-Cauchy.
A subset $E$ of $\textbf{R}$ is called downward compact if any sequence of points in $E$ has a downward half quasi-Cauchy subsequence.
A subset of $\textbf{R}$ is downward compact if and only if it is bounded above.
A subset $E$ of $\textbf{R}$ is called down half compact if any sequence of points in $E$ has an down half Cauchy subsequence (\cite{CakalliHalfquasi-Cauchysequences}).
A subset of $\textbf{R}$ is down half compact if and only if it is downward compact .
A subset of $\textbf{R}$ is bounded if and only if it is both upward and downward compact.

\maketitle

\section{Statistically upward and downward quasi-Cauchy sequences}

Weakening the condition on the definition of a statistically quasi-Cauchy sequence, omitting the absolute value symbol, i.e. if we replace $|x_{k}-x_{k+1}|\geq \varepsilon$ by $x_{k}-x_{k+1} \geq \varepsilon$ in the definition of a statistically quasi-Cauchy sequence, we have the following definition.

\begin{Def} A sequence $(x_{n})$ of points in $\textbf{R}$ is called statistically upward half quasi-Cauchy if
\[
\lim_{n\rightarrow\infty}\frac{1}{n}|\{k\leq n: x_{k}-x_{k+1}\geq \varepsilon\}|=0
\]
for every $\varepsilon>0$.
\end{Def}

Now we introduce a definition of statistically upward compactness of a subset of $\textbf{R}$, by using the main idea in the definition of sequential compactness.

\begin{Def}
A subset $E$ of $\textbf{R}$ is called statistically upward compact if any sequence of points in $E$ has an statistically upward half quasi-Cauchy subsequence.
\end{Def}

First, we note that any finite subset of $\textbf{R}$ is statistically upward compact, union of two statistically upward compact subsets of $\textbf{R}$ is statistically upward compact and intersection of any statistically upward compact subsets of $\textbf{R}$ is statistically upward compact. Furthermore any subset of a statistically upward compact set is statistically upward compact and any bounded subset of $\textbf{R}$ is statistically upward compact. Any compact subset of $\textbf{R}$ is also statistically upward compact. We note that any slowly oscillating compact subset of $\textbf{R}$ is statistically upward compact (see [7] for the definition of slowly oscillating compactness). These observations suggest to us the following:

\begin{Thm} \label{Theohalfstatisticallyupwardcompactiffboundedbelow}

A subset of $\textbf{R}$ is statistically upward compact if and only if it is bounded below.

\end{Thm}

\begin{pf}
Let $E$ be a subset of $\textbf{R}$. If $E$ is bounded below, then it is clear that $E$ is statistically upward compact.
Now, suppose that $E$ is not bounded below. Pick an element $x_{1}$ of $E$. Then we can choose an element $x_{2}$ of $E$ such that $x_{2}<-1+x_{1}$. Similarly we can choose an element $x_{3}$ of $E$ such that $x_{3}<-1+x_{2}$. We can inductively choose  $k_{n}$ and $k_{n+1}$ satisfying
$x_{n+1}<-1+x_{n}$ for each positive integer $n$. Then the sequence $(x_{n})$ does not have any statistically upward half quasi-Cauchy subsequence. Thus $E$ is not statistically upward  compact. This contradiction completes the proof.
\end{pf}

\begin{Def} A sequence $(x_{n})$ of points in $\textbf{R}$ is called statistically downward half quasi-Cauchy if
\[
\lim_{n\rightarrow\infty}\frac{1}{n}|\{k\leq n: x_{k+1}-x_{k}\geq \varepsilon\}|=0
\]
for every $\varepsilon>0$.
\end{Def}
Throughout the paper $\Delta S^{+}$ and $\Delta S^{-}$  will denote the set of all statistically upward half quasi-Cauchy sequences, and the set of all statistically downward half quasi-Cauchy sequences, of points in $\textbf{R}$, respectively. Any statistically quasi-Cauchy sequence is statistically upward half quasi-Cauchy, statistically downward half quasi-Cauchy, so any slowly oscillating sequence is statistically upward and downward half quasi-Cauchy, so any Cauchy sequence is, so any convergent sequence is. Any upward half Cauchy sequence is statistically upward half quasi-Cauchy.

Now we introduce a definition of statistically downward compactness of a subset of $\textbf{R}$.

\begin{Def}
A subset $E$ of $\textbf{R}$ is called statistically downward compact if any sequence of points in $E$ has a statistically downward half quasi-Cauchy subsequence.
\end{Def}

First, we note that any finite subset of $\textbf{R}$ is statistically downward compact, union of two statistically downward compact subsets of $\textbf{R}$ is statistically downward compact and intersection of any statistically downward compact subsets of $\textbf{R}$ is statistically downward compact. Furthermore any subset of a statistically downward compact set is statistically downward compact and any bounded subset of $\textbf{R}$ is statistically downward compact. Any compact subset of $\textbf{R}$ is also statistically downward compact. We note that any slowly oscillating compact subset of $\textbf{R}$ is statistically downward compact (see [7] for the definition of slowly oscillating compactness).

\begin{Thm} \label{Theohalfupwardcompactiffboundedbelow}

A subset of $\textbf{R}$ is statistically downward compact if and only if it is bounded above.

\end{Thm}

\begin{pf}
Let $E$ be a subset of $\textbf{R}$. If $E$ is bounded above, then it is clear that $E$ is statistically downward compact.
Now, suppose that $E$ is not bounded above. Pick an element $x_{1}$ of $E$. Then we can choose an element $x_{2}$ of $E$ such that $x_{2}>1+x_{1}$. Similarly we can choose an element $x_{3}$ of $E$ such that $x_{3}>1+x_{2}$. We can inductively choose  $k_{n}$ and $k_{n+1}$ satisfying
$x_{n+1}>1+x_{n}$ for each positive integer $n$. Then the sequence $(x_{n})$ does not have any statistically downward half quasi-Cauchy subsequence. Thus $E$ is not statistically downward compact. This is a contradiction.

\end{pf}

A subset $E$ of $\textbf{R}$ is called down half compact if any sequence of points in $E$ has an down half Cauchy subsequence (\cite{CakalliHalfquasi-Cauchysequences}).

\begin{Cor} \label{Corhalfdownhalfcompactiffstatisticallydownwardcompact}
A subset of $\textbf{R}$ is down half compact if and only if it is statistically downward compact .
\end{Cor}

\begin{Cor} \label{CorboundednesscoincideswithstatisupwardquasiandhalfdownwardquasiCauchycompactness}

A subset of $\textbf{R}$ is bounded if and only if it is both statistically upward and statistically downward compact.

\end{Cor}

\begin{Cor} \label{Corstatisticalwardcompactnesscoincideswithstatiswardandstatisdownwardcompactness}

A subset of $\textbf{R}$ is statistically ward compact if and only if it is both statistically upward and statistically downward compact.

\end{Cor}

\begin{Cor} \label{CorNthetawardcompactnesscoincideswistatisupwardcompactnessandstatisdownwardcompactness}

A subset of $\textbf{R}$ is $N_{\theta}$-ward compact if and only if it is both statistically upward and statistically downward compact.

\end{Cor}

\begin{pf} The proof follows from \ref{CorboundednesscoincideswithstatisupwardquasiandhalfdownwardquasiCauchycompactness}, and \cite [Theorem 3.3]{CakalliNthetawardcontinuity}.
\end{pf}

\begin{Def} A sequence $(x_{n})$ of points in $\textbf{R}$ is called half statistically quasi-Cauchy if the sequence is either statistically  upward, or statistically downward half quasi-Cauchy, or both.
\end{Def}

We note that a sequence $(x_{n})$ of points in $\textbf{R}$ is statistically quasi-Cauchy if and only if it is both statistically upward half quasi-Cauchy and statistically downward half quasi-Cauchy.

\section{Statistically Upward and Downward continuities}

A real valued function $f$ defined on a subset of $\textbf{R}$ is statistically continuous if and only if, for each point $\ell$ in the domain, $st-\lim_{n\rightarrow\infty}f(x_{n})=f(\ell)$ whenever $st-\lim_{n\rightarrow\infty}x_{n}=\ell$. This is equivalent to the statement that $(f(x_{n}))$ is a statistically convergent sequence whenever $(x_{n})$ is. This is also equivalent to the statement that $(f(x_{n}))$ is a Cauchy sequence whenever $(x_{n})$ is Cauchy. These well known results for statistical continuity for real functions in terms of sequences might suggest to us giving a new type continuity, namely, statistically upward continuity:

Now we give the concept of statistically upward continuity of a function defined on a subset of $\textbf{R}$ to $\textbf{R}$.

\begin{Def}
A function $f$ is called statistically upward continuous on $E$ if the sequence $(f(x_{n}))$ is statistically upward half quasi-Cauchy whenever $\textbf{x}=(x_{n})$ is a statistically upward half quasi-Cauchy sequence of points in $E$.
\end{Def}

We see that sum of two statistically upward continuous functions is statistically upward  continuous and composite of two statistically upward continuous functions is statistically upward  continuous.

In connection with statistically upward half quasi-Cauchy sequences and convergent sequences the problem arises to investigate the following types of  "continuity" of functions on $\textbf{R}$.

\begin{description}
\item[($\delta S^{+} $)] $(x_{n}) \in {\Delta S^{+}} \Rightarrow (f(x_{n})) \in {\Delta S^{+}}$
\item[($\delta S^{+} c$)] $(x_{n}) \in {\Delta S^{+}} \Rightarrow (f(x_{n})) \in {c}$
\item[$(c)$] $(x_{n}) \in {c} \Rightarrow (f(x_{n})) \in {c}$
\item[$(c\delta S^{+})$] $(x_{n}) \in {c} \Rightarrow (f(x_{n})) \in {\Delta S^{+}}$
\item[$(st)$] $(x_{n}) \in {S} \Rightarrow (f(x_{n})) \in {S}$
\end{description}

We see that $(\delta S^{+})$ is statistically upward continuity of $f$, and $(st)$ is the statistical continuity. It is easy to see that $(\delta S^{+} c)$ implies $(\delta S^{+})$; $(\delta S^{+})$ does not imply $(\delta S^{+} c)$; $(\delta S^{+})$ implies $(c\delta S^{+})$; $(c\delta S^{+})$ does not imply $(\delta S^{+})$; $(\delta S^{+} c)$ implies $(c)$, and $(c)$ does not imply $(\delta S^{+} c)$; and $(c)$ is equivalent to $(c\delta S^{+})$.

Now we prove that ($\delta S^{+} $) implies $(c)$ in the following:

\begin{Thm} \label{TheoStatisticallyupwardcontinuityimpliesstatisticalcontinuity} If $f$ is statistically upward continuous on a subset $E$ of $\textbf{R}$, then it is statistically continuous on $E$.
\end{Thm}
\begin{pf}

Let $(x_{n})$ be any statistically convergent sequence with $st-\lim_{k\rightarrow\infty}x_{k}=\ell$. Then $$(x_{1}, \ell , x_{1}, \ell , x_{2}, \ell, x_{2}, \ell,..., x_{n}, \ell, x_{n}, \ell,...)$$ is also statistically convergent to $\ell$. Thus it is statistically upward half quasi-Cauchy. Hence
$$(f(x_{1}), f(\ell), f(x_{1}),f(\ell), f(x_{2}), f(\ell), f(x_{2}), f(\ell),...,f(x_{n}),f(\ell),f(x_{n}),f(\ell),...)$$
is statistically upward half quasi-Cauchy. It follows from this that
\[
\lim_{n\rightarrow\infty}\frac{1}{n}|\{k\leq n: f(x_{k})-f(\ell)\geq \varepsilon \; \; \; \normalfont{and} \; \; \; f(x_{k+1})-f(\ell)\geq \varepsilon\}|=0
\]
for every $\varepsilon>0$
which implies that
\[
\lim_{n\rightarrow\infty}\frac{1}{n}|\{k\leq n: |f(x_{k})-f(\ell)| \geq \varepsilon \}|=0
\]
for every $\varepsilon>0$. This completes the proof of the theorem.
\end{pf}

Now we state the following result related to ordinary continuity and statistically upward continuity.

\begin{Cor} If $f$ is statistically upward continuous, then it is ordinary continuous.
\end{Cor}

\begin{Cor} If $f$ is statistically upward continuous, then it is lacunary statistically continuous.
\end{Cor}

\begin{Cor} If $f$ is statistically upward continuous, then it is $N_{\theta}$-continuous.

\end{Cor}

We have the following result for general sequential methods.

\begin{Cor} If $f$ is statistically upward continuous, then it is $G$-continuous for any regular subsequential method $G$.

\end{Cor}

\begin{Cor} If $f$ is statistically upward continuous, then it is $I$-continuous for any non-trivial admissible ideal $I$ of $\textbf{N}$.

\end{Cor}

\begin{pf}
The proof follows from \ref{TheoStatisticallyupwardcontinuityimpliesstatisticalcontinuity}, and \cite[Theorem 4]{CakalliandHazarikaIdealquasi-Cauchysequences}.
\end{pf}

\begin{Thm} \label{TheoStatisticallyupwardcontinousimageofstatisticallyupwardcompactsubsetis} Statistically upward continuous image of any statistically upward compact subset of $\textbf{R}$ is statistically upward compact.
\end{Thm}
\begin{pf}
Write $y_{n}=f(x_{n})$ where $x_{n}\in {E}$ for each $n \in{\textbf{N}}$. Statistically upward compactness of $E$ implies that there is a statistically upward half quasi-Cauchy subsequence of the sequence of $\textbf{x}$. Write $(t_{k})=f(\textbf{z})=(f(z_{k}))$. $(t_{k})$ is a statistically upward half quasi-Cauchy subsequence of the sequence $f(\textbf{x})$. This completes the proof of the theorem.
\end{pf}

\begin{Cor} Statistically upward continuous image of any $N_{\theta}$-sequentially compact subset of $\textbf{R}$ is $N_{\theta}$-sequentially compact.
\end{Cor}

\begin{Cor} Statistically upward continuous image of any compact subset of $\textbf{R}$ is compact.
\end{Cor}

\begin{Thm} Let $E$ be a statistically upward compact subset $E$ of \textbf{R} and let $f:E\longrightarrow$ $\textbf{R}$ be a  statistically upward continuous function on $E$. Then $f$ is uniformly continuous on $E$.
\end{Thm}

\begin{pf}
Suppose that $f$ is not uniformly continuous on $E$ so that there exists an  $\varepsilon_{0} > 0$ such that for any $\delta >0$ $x, y \in{E}$ with $|x-y|<\delta$ but $|f(x)-f(y)| \geq \varepsilon_{0}$. For each positive integer $n$, fix $|x_{n}-y_{n}|<\frac{1}{n}$, and $|f(x_{n})-f(y_{n})|\geq \varepsilon_{0}$. Since $E$ is statistically upward compact, there exists a statistical upward half quasi-Cauchy subsequence $(x_{n_{k}})$ of the sequence $(x_{n})$. It is clear that the corresponding subsequence $(y_{n_{k}})$ of the sequence $(y_{n})$ is also statistically upward half quasi-Cauchy, since $(y_{n_{k}}-y_{n_{k+1}})$ is a sum of three statistically upward half quasi-Cauchy sequences, i.e. $$y_{n_{k}}-y_{n_{k+1}}=(y_{n_{k}}-x_{n_{k}})+(x_{n_{k}}-x_{n_{k+1}})+(x_{n_{k+1}}-y_{n_{k+1}}).$$ Then the sequence
$$(x_{n_{1}}, y_{n_{1}}, x_{n_{2}}, y_{n_{2}}, x_{n_{3}}, y_{n_{3}},..., x_{n_{k}}, y_{n_{k}},...)$$
is statistically upward quasi-Cauchy since the sequence $(x_{n_{k}}-y_{n_{k+1}})$ is statistically upward quasi-Cauchy sequence which follows from the equality $$x_{n_{k}}-y_{n_{k+1}}=x_{n_{k}}-x_{n_{k+1}}+x_{n_{k+1}}-y_{n_{k+1}}.$$ But the transformed sequence
$$(f(x_{n_{1}}), f(y_{n_{1}}), f(x_{n_{2}}), f(y_{n_{2}}), f(x_{n_{3}}), f(y_{n_{3}}),..., f(x_{n_{k}}), f(y_{n_{k}}),...)$$
is not statistically upward half quasi-Cauchy. Thus $f$ does not preserve statistically upward half quasi-Cauchy sequences. This contradiction completes the proof of the theorem.
\end{pf}

It is a well known result that uniform limit of a sequence of continuous functions is continuous. This is also true in case statistically  upward continuity, i.e. uniform limit of a sequence of statistically upward continuous functions is statistically upward continuous.
\begin{Thm} If $(f_{n})$ is a sequence of statistically upward continuous functions defined on a subset $E$ of $\textbf{R}$ and $(f_{n})$ is uniformly convergent to a function $f$, then $f$ is statistically upward continuous on $E$.
\end{Thm}
\begin{pf}
Let $\varepsilon$ be a positive real number and $(x_{k})$ be any statistically upward half quasi-Cauchy sequence of points in $E$. By uniform convergence of $(f_{n})$ there exists a positive integer $N$ such that $|f_{n}(x)-f(x)|<\frac{\varepsilon}{3}$ for all $x \in {E}$ whenever $n\geq N$. As $f_{N}$ is statistically upward continuous on $E$, we have $$\lim_{n\rightarrow\infty} \frac{1}{n} |\{k \leq {n} : f_{N}(x_{k})-f_{N}(x_{k+1})\geq \frac{\varepsilon}{3}\}|=0.$$ On the other hand we have\\
$\{ k \leq {n} : f(x_{k})-f(x_{k+1})\geq \varepsilon\} \subset {\{k \leq {n} : f(x_{k})-f_{N}(x_{k}) \geq \frac{\varepsilon}{3}\}}$
 \; \; \; \; \; \; \; \; \; \; \; $\cup \{k \leq {n} :f_{N}(x_{k})-f_{N}(x_{k+1}) \geq \frac{\varepsilon}{3}\} \cup \{k \leq {n} :f_{N}(x_{k+1})-f(x_{k+1}) \geq \frac{\varepsilon}{3}\}$ \\
 and so\\
 $\{ k \leq {n} : f(x_{k})-f(x_{k+1})\geq \varepsilon\} \subset {\{k \leq {n} : |f(x_{k})-f_{N}(x_{k})|\geq \frac{\varepsilon}{3}\}}$
 \; \; \; \; \; \; \; \; \; \; \; $\cup \{k \leq {n} :f_{N}(x_{k})-f_{N}(x_{k+1}) \geq \frac{\varepsilon}{3}\} \cup \{k \leq {n} :|f_{N}(x_{k+1})-f(x_{k+1})| \geq \frac{\varepsilon}{3}\}$. \\
Now it follows from this inclusion that \\
$ \lim_{n\rightarrow\infty}\frac{1}{n} | \{ k \leq {n} : f(x_{k})-f(x_{k+1}) \geq \varepsilon\}|$
 $$\leq   \lim_{n\rightarrow\infty}\frac{1}{n} |{\{k \leq {n} : |f(x_{k})-f_{N}(x_{k})| \geq \frac{\varepsilon}{3}\}| + \lim_{n\rightarrow\infty}\frac{1}{n} |\{k \leq {n} : f_{N}(x_{k})-f_{N}(x_{k+1}) \geq \frac{\varepsilon}{3}}\}|$$   $$ + \lim_{n\rightarrow\infty}\frac{1}{n} |\{k \leq {n} :|f_{N}(x_{k+1})-f(x_{k+1})| \geq \frac{\varepsilon}{3}\}|=0+0+0=0.$$
This completes the proof of the theorem.
\end{pf}

Now we give the definition of statistically downward continuity of a real function in the following:

\begin{Def}
A function $f$ is called statistically downward continuous on $E$ if the sequence $(f(x_{n}))$ is statistically downward half Cauchy whenever $\textbf{x}=(x_{n})$ is a statistically downward half Cauchy sequence of points in $E$.
\end{Def}

We note that sum of two statistically downward continuous functions is statistically downward continuous and composite of two downward continuous functions is downward continuous.

Now we investigate the following types of  "continuity" of functions on $\textbf{R}$.

\begin{description}
\item[($\delta S^{-}$)] $(x_{n}) \in {\Delta S^{-}} \Rightarrow (f(x_{n})) \in {\Delta S^{-}}$
\item[($\delta S^{-}c$)] $(x_{n}) \in {\Delta S^{-}} \Rightarrow (f(x_{n})) \in {c}$
\item[$(c\delta S^{-})$] $(x_{n}) \in {c} \Rightarrow (f(x_{n})) \in {\Delta S^{-}}$
\end{description}

We see that $(\delta S^{-})$ is downward statistical continuity of $f$. It is easy to see that $(\delta S^{-}c)$ implies $(\delta S^{-})$; $(\delta S^{-})$ does not imply $(\delta S^{-}c)$; $(\delta S^{-})$ implies $(c\delta S^{-})$; $(c\delta S^{-})$ does not imply $(\delta S^{-})$; $(\delta S^{-}c)$ implies $(c)$; $(c)$ does not imply $(\delta S^{-}c)$, and $(c)$ is equivalent to $(c\delta S^{-})$.

Now we give the implication $(\delta S^{-})$ implies $(st)$, i.e. any statistically downward continuous function is statistically continuous.

\begin{Thm} \label{TheoStatisticallydownwardcontinuityimpliesstatisticalcontinuity}  If $f$ is statistically downward continuous on a subset $E$ of $\textbf{R}$, then it is statistically continuous on $E$.
\end{Thm}
\begin{pf}
Let $(x_{n})$ be any statistically convergent sequence with $st-\lim_{k\rightarrow\infty}x_{k}=\ell$. Then $$(x_{1}, \ell , x_{1}, \ell , x_{2}, \ell, x_{2}, \ell,..., x_{n}, \ell, x_{n}, \ell,...)$$ is also statistically convergent to $\ell$. Thus it is statistically downward half quasi-Cauchy. Hence $$(f(x_{1}), f(\ell), f(x_{1}),f(\ell), f(x_{2}), f(\ell), f(x_{2}), f(\ell),...,f(x_{n}), f(\ell), f(x_{n}), f(\ell),...)$$ is statistically downward half quasi-Cauchy. It follows from this that
\[
\lim_{n\rightarrow\infty}\frac{1}{n}|\{k\leq n: f(x_{k+1})-f(\ell)\geq \varepsilon \; \; \; \normalfont{and} \; \; \; f(x_{k})-f(\ell)\geq \varepsilon\}|=0
\]
for every $\varepsilon>0$
which implies that
\[
\lim_{n\rightarrow\infty}\frac{1}{n}|\{k\leq n: |f(x_{k})-f(\ell)| \geq \varepsilon \}|=0
\]
for every $\varepsilon>0$. This completes the proof of the theorem.

\end{pf}

\begin{Cor} If $f$ is statistically downward continuous, then it is continuous.
\end{Cor}

\begin{Cor} If $f$ is statistically downward continuous, then it is lacunary statistically continuous.
\end{Cor}

\begin{Cor} If $f$ is statistically downward continuous, then it is $N_{\theta}$-continuous.

\end{Cor}

\begin{Cor} If $f$ is statistically downward continuous, then it is $G$-continuous for any regular subsequential method $G$.

\end{Cor}

\begin{Cor} If $f$ is statistically downward continuous, then it is $I$-continuous for any non-trivial admissible ideal $I$ of $\textbf{N}$.

\end{Cor}

\begin{pf}
The proof follows from \ref{TheoStatisticallydownwardcontinuityimpliesstatisticalcontinuity}, and \cite[Theorem 4]{CakalliandHazarikaIdealquasi-Cauchysequences}.
\end{pf}

\begin{Thm} \label{ThehalfquasiCauchycontinousimageofhalfquasi-Cauchycompactsubsetis} Statistically downward continuous image of any statistically downward half compact subset of $\textbf{R}$ is statistically downward half compact.
\end{Thm}
\begin{pf}
Write $y_{n}=f(x_{n})$ where $x_{n}\in {E}$ for each $n \in{\textbf{N}}$. Statistically downward half compactness of $E$ implies that there is a statistically downward half quasi-Cauchy subsequence $\textbf{z}=(z_{k})=(x_{n_{k}})$ of $\textbf{x}$. Write $(t_{k})=f(\textbf{z})=(f(z_{k}))$. Then $(t_{k})$ is a statistically downward half quasi-Cauchy subsequence of the sequence $f(\textbf{x})$. This completes the proof of the theorem.
\end{pf}

\begin{Cor} Statistically downward continuous image of any $N_{\theta}$-sequentially compact subset of $\textbf{R}$ is $N_{\theta}$-sequentially compact.
\end{Cor}

\begin{Cor} Statistically downward image of any compact subset of $\textbf{R}$ is compact.
\end{Cor}

\begin{Thm} Let $E$ be a statistically downward compact subset $E$ of \textbf{R} and let $f:E\longrightarrow$ $\textbf{R}$ be a  statistically downward continuous function on $E$. Then $f$ is uniformly continuous on $E$.
\end{Thm}

\begin{pf}
Suppose that $f$ is not uniformly continuous on $E$ so that there exists an  $\varepsilon_{0} > 0$ such that for any $\delta >0$ $x, y \in{E}$ with $|x-y|<\delta$ but $|f(x)-f(y)| \geq \varepsilon_{0}$. For each positive integer $n$, fix $|x_{n}-y_{n}|<\frac{1}{n}$, and $|f(x_{n})-f(y_{n})|\geq \varepsilon_{0}$. Since $E$ is statistically downward compact, there exists a statistical downward half quasi-Cauchy subsequence $(x_{n_{k}})$ of the sequence $(x_{n})$. It is clear that the corresponding subsequence $(y_{n_{k}})$ of the sequence $(y_{n})$ is also statistically downward half quasi-Cauchy, since $(y_{n_{k+1}}-y_{n_{k}})$ is a sum of three statistical downward half quasi-Cauchy sequences, i.e. $$y_{n_{k+1}}-y_{n_{k}}=(y_{n_{k+1}}-x_{n_{k+1}})+(x_{n_{k+1}}-x_{n_{k}})+(x_{n_{k}}-y_{n_{k}}).$$ Then the sequence
$$(y_{n_{1}}, x_{n_{1}},  y_{n_{2}}, x_{n_{2}}, y_{n_{3}}, x_{n_{3}}, ..., y_{n_{k}}, x_{n_{k}}, ...)$$
is statistical downward half quasi-Cauchy since the sequence $(y_{n_{k+1}}-x_{n_{k}})$ is statistically downward half quasi-Cauchy sequence which follows from the equality $$y_{n_{k+1}}-x_{n_{k}}=y_{n_{k+1}}-y_{n_{k}}+y_{n_{k}}-x_{n_{k}}$$. But the transformed sequence
$$(f(y_{n_{1}}), f(x_{n_{1}}), f(y_{n_{2}}), f(x_{n_{2}}), f(y_{n_{3}}), f(x_{n_{3}}),..., f(y_{n_{k}}), f(x_{n_{k}}),...)$$
is not statistically downward half quasi-Cauchy. Thus $f$ does not preserve statistically downward half quasi-Cauchy sequences. This contradiction completes the proof of the theorem.
\end{pf}

Now we have the following result related to uniform convergence, namely, uniform limit of a sequence of statistically downward continuous functions is again statistically downward continuous.

\begin{Thm} If $(f_{n})$ is a sequence of statistically downward continuous functions defined on a subset $E$ of $\textbf{R}$ and $(f_{n})$ is uniformly convergent to a function $f$, then $f$ is statistically downward continuous on $E$.
\end{Thm}
\begin{pf}
 Let $\varepsilon$ be a positive real number and $(x_{k})$ be any statistically downward half quasi-Cauchy sequence of points in $E$. By uniform convergence of $(f_{n})$ there exists a positive integer $N$ such that $|f_{n}(x)-f(x)|<\frac{\varepsilon}{3}$ for all $x \in {E}$ whenever $n\geq N$. As $f_{N}$ is statistically downward continuous on $E$, we have $$\lim_{n\rightarrow\infty} \frac{1}{n} |\{k \leq {n} : f_{N}(x_{k+1})-f_{N}(x_{k}) \geq \frac{\varepsilon}{3}\}|=0.$$ On the other hand we have
 \\
$\{ k \leq {n} :f(x_{k+1})-f(x_{k}) \geq \varepsilon\} = {\{k \leq {n} : f(x_{k+1})-f_{N}(x_{k+1})|\geq \frac{\varepsilon}{3}\}}$
 \; \; \; \; \; \; \; \; \; \; \; $\cup \{k \leq {n} :f_{N}(x_{k+1})-f_{N}(x_{k}) \geq \frac{\varepsilon}{3}\} \cup \{k \leq {n} : f_{N}(x_{k})-f(x_{k}) \geq \frac{\varepsilon}{3}\}$
 \\so we have\\
$\{ k \leq {n} :f(x_{k+1})-f(x_{k}) \geq \varepsilon\} \subset {\{k \leq {n} : |f(x_{k+1})-f_{N}(x_{k+1})|\geq \frac{\varepsilon}{3}\}}$
 \; \; \; \; \; \; \; \; \; \; \; $\cup \{k \leq {n} :f_{N}(x_{k+1})-f_{N}(x_{k}) \geq \frac{\varepsilon}{3}\} \cup \{k \leq {n} :|f_{N}(x_{k})-f(x_{k})| \geq \frac{\varepsilon}{3}\}$ \\
Now it follows from this inclusion that \\
$ \lim_{n\rightarrow\infty}\frac{1}{n} | \{ k \leq {n} :f(x_{k+1})-f(x_{k}) \geq \varepsilon\}|$
 $$\leq   \lim_{n\rightarrow\infty}\frac{1}{n} |{\{k \leq {n} : |f(x_{k+1})-f_{N}(x_{k+1})| \geq \frac{\varepsilon}{3}\}| + \lim_{n\rightarrow\infty}\frac{1}{n} |\{k \leq {n} : f_{N}(x_{k+1})-f_{N}(x_{k}) \geq \frac{\varepsilon}{3}}\}|$$   $$ + \lim_{n\rightarrow\infty}\frac{1}{n} |\{k \leq {n} :|f_{N}(x_{k})-f(x_{k})| \geq \frac{\varepsilon}{3}\}|=0+0+0=0.$$
This completes the proof of the theorem.
\end{pf}

\maketitle

\section{Conclusion}

In this paper, we introduce and investigate statistically upward and statistically downward continuities of real functions, statistically upward and statistically downward compactnesses of a subset of $\textbf{R}$. We prove results related to these kinds of continuities, these kind of compactnesses, compactness,  and some other kinds of continuities and compactness; namely slowly oscillating continuity, slowly oscillating compactness, statistical compactness, lacunary statistical compactness, ordinary compactness, ordinary continuity, and uniform continuity. It turns out that not only the set of upward continuous functions on upward compact subsets but also the set of downward continuous functions on downward compact subsets coincides with the set of uniformly continuous functions, uniform limit of a sequence of upward continuous functions is upward continuous, and uniform limit of a sequence of downward continuous functions is downward continuous.
We suggest to investigate upward and downward half quasi-Cauchy sequences of fuzzy points or soft points (see \cite{CakalliandPratul} for the definitions and  related concepts in fuzzy setting, and see \cite{CagmanandKaratasandEnginoglu[19]Softtopology}, and \cite{ArasandSonmezandCakalliOnSoftMappings} related concepts in soft setting), and upward and downward half quasi-Cauchy sequences in cone metric spaces (\cite{SonmezandCakalliConenormedspacesandweightedmeans}, \cite{CakalliandSonmezandGenc} and \cite{PalandSavasandCakalliIconvergenceonconemetricspaces}). However due to the change in settings, the definitions and methods of proofs will not always be analogous to those of the present work. For some other further studies we suggest to investigate upward and downward half quasi-Cauchy sequences of double sequences (see for example \cite{MursaleenandMohiuddineBanachlimitandsomenewspacesofdoublesequences}, and \cite{PattersonandSavasAsymptoticEquivalenceofDoubleSequences} for the definitions and related concepts in the double case).

\end{document}